\documentclass{amsart}
\usepackage[english]{babel}
\usepackage[latin1]{inputenc}
\usepackage{mathrsfs}
\usepackage{amsmath}
\usepackage{amsthm}
\usepackage{amssymb}
\usepackage{indentfirst}

\usepackage[dvips]{graphicx}

\usepackage{ifthen}
\newcommand{\cit}[1]{{\rm \textbf{#1}}}

\newcommand{\Ref}[2]{\cit{%
\ifthenelse{\equal{#1}{thm}}{Theorem}{}%
\ifthenelse{\equal{#1}{prop}}{Proposition}{}%
\ifthenelse{\equal{#1}{lem}}{Lemma}{}%
\ifthenelse{\equal{#1}{cor}}{Corollary}{}%
\ifthenelse{\equal{#1}{defn}}{Definition}{}%
\ifthenelse{\equal{#1}{oss}}{Remark}{}%
\ifthenelse{\equal{#1}{sec}}{Section}{}%
\ifthenelse{\equal{#1}{ex}}{Example}{}%
\ifthenelse{\equal{#1}{conj}}{Conjecture}{}%
\ifthenelse{\equal{#1}{ssec}}{Subsection}{}%
\ifthenelse{\equal{#1}{tab}}{Table}{}%
\ifthenelse{\equal{#1}{cla}}{Claim}{}%
\  \ref{#1:#2}%
}}

\newtheorem{prop}{Proposition}[section]
\newtheorem{thm}[prop]{Theorem}
\newtheorem{lem}[prop]{Lemma} 
\newtheorem{cor}[prop]{Corollary}

\theoremstyle{definition}
\newtheorem{defn}[prop]{Definition}
\newtheorem{ex}[prop]{Example}

\theoremstyle{remark}
\newtheorem{oss}[prop]{Remark}

\numberwithin{equation}{section}

\newcommand{\hk}{hyperk\"{a}hler }

\newcommand{\kahl}{K\"{a}hler }
\newcommand{\ktiposp}{$K3^{[2]}$ type }
\newcommand{\ktipo}{$K3^{[2]}$ type}
\newcommand{\kntipo}{$K3^{[n]}$ type}
\newcommand{\kntiposp}{$K3^{[n]}$ type }
\newcommand{\ie}{i.~e.~}
\newcommand{\muraja}{wall divisor}
\newcommand{\murajasp}{wall divisor }
\newcommand{\muraje}{wall divisors}
\newcommand{\murajesp}{wall divisors }

\begin{document}

\title{A note on the K\"{a}hler and Mori cones of hyperk\"{a}hler manifolds}

\author{Giovanni Mongardi}
\address{Mathematisches Institut der Universit\"{a}t Bonn. Endenicher Allee, 60}
\curraddr{Universit\`{a} degli studi di Milano, Dipartimento di Matematica, via Cesare Saldini, 50}
\email{giovanni.mongardi@unimi.it}
\thanks{Supported by the SFB/TR 45 ``Periods, Moduli Spaces and Arithmetic of Algebraic Varieties''}


\date{}


\begin{abstract}
In  the present paper we prove that, on a \hk manifold, walls of the \kahl cone and extremal rays of the Mori cone are determined by all divisors satisfying certain numerical conditions.
\end{abstract}

\maketitle

\section*{Introduction}

The structure of the Mori cone and the \kahl cones of \hk manifolds has attracted some attention in recent years. For $K3$ surfaces, the Mori cone is generated by curves of self intersection $\geq-2$ and the \kahl cone is just the dual cone. Huybrechts proved some pioneering results for general \hk manifolds in \cite{huy_ann} and more was done by Markman (cf. \cite{mar2}). In the present paper we introduce a class of divisors, \muraje, preserved under deformations where they stay of $(1,1)$ type. This class contains divisors associated to extremal rays of the Mori cone, cf. \Ref{prop}{mori}. This class of divisors is preserved under smooth deformations, which is the main result of the paper (cf. \Ref{thm}{muri_kahl}) and gives a wall and chamber decomposition of the positive cone. In particular, once a way to determine all possible orbits of extremal rays is given, one of the chambers in this decomposition is the \kahl cone.\\  In the case of manifolds of \kntiposp it had been conjectured by Hassett and Tschinkel that the Mori cone is given by curves with self intersection greater than $-\frac{n+3}{2}$. Using minimal model program, they proved in \cite{ht2} one direction of the conjecture in dimension four. Although there are now counterexamples to the general form of this conjecture (see \cite[Proposition 9.3]{bm2}), it is still expected that some numerical characterization of the Mori cone can be given. For the closure of the birational \kahl cone this numerical characterization has been given by Markman \cite{mar2}, we summarize it in \Ref{thm}{mark_muri}. Bayer and Macr\`{i} \cite{bm} give a numerical characterization of the \kahl and Mori cones for manifolds which are given by moduli spaces of sheaves on projective $K3$ surfaces, we summarize their results in \Ref{thm}{eman_muri}.\\ As an application of our class of divisors, we use the latter result to give a numerical characterization for all manifolds of \kntipo, including non algebraic ones. We analyze several divisors known to be \murajesp and we obtain a full numerical classification of the \kahl cone for $n=2,3$ and $4$. It is of course possible to extend these results to arbitrary $n$ using the known classification on Moduli spaces. For all manifolds of \kntipo, we prove one direction of the Hassett and Tschinkel conjecture, namely that extremal rays of the Mori cone have square greater than $-\frac{n+3}{2}$, thus generalizing the proof of Bayer and Macr\`{i} to arbitrary deformations.\\ 
A nice geometrical interpretation of some walls of the \kahl cone can be found in \cite{ht4}, where Hasset and Tschinkel analyzed contractions of some extremal rays of the Mori cone. It should be pointed out that, by a result of Markman \cite{mar2}, most of the divisors dual to extremal rays are not irreducible and, in the non algebraic case, they are also not effective. This gives an example of analytic manifolds with no effective divisors but with rational curves.     

\subsection*{Acknowledgements}
I would like to thank Prof.$\;$D. Huybrechts for his advice, Prof.$\;$E. Macr\`{i} and Prof.$\;$A. Bayer for their comments. They also informed me that a different proof of an equivalent result in the case of projective manifolds of \kntiposp is proven in \cite{bht} by Prof.$\;$Bayer, Prof.$\;$Hassett and Prof.$\;$Tschinkel and I would once more like to thank Prof.$\;$Bayer for a preliminary version of it. I am grateful to Prof. E. Markman for a detailed discussion on isometries preserving wall divisors and for his many comments. I would also like to thank M. Lelli-Chiesa for reading a preliminary version of this paper and U. Greiner for an insightful question. I am also grateful to the referee for helpful comments. This paper has been revised under the support of FIRB 2012 ``Spazi di Moduli e applicazioni'' and of the Max Planck Institute for Mathematics in Bonn.

\subsection*{Notations}
Let $M$ be a lattice of signature $(3,r)$. Let $X$ be a \hk manifold and let $f:\,H^2(X,\mathbb{Z})\,\rightarrow\,M$ be a marking. The moduli space of marked \hk manifolds is denoted by $\mathcal{M}_M$. The connected component containing $(X,f)$ is denoted by $\mathcal{M}_M^0$.
The period domain $\Omega_M$ is defined in two equivalent ways: as the Grassmannian of positive oriented 2 planes inside $M\otimes\mathbb{R}$ or the open set $\{\omega\in\mathbb{P}(M\otimes\mathbb{C}),\,\omega^2=0,\,(\omega,\overline{\omega})>0\}$. In the first case the period map $\mathcal{P}$ sends $(X,f)$ to the positive plane $\langle\,Re(f(\sigma_X)), Im(f(\sigma_X))\,\rangle$, while in the second the image is $[f(\sigma_X)]$.\\
Let $g$ be a \hk metric on a manifold $X$ with \kahl class $\omega$. The \emph{twistor family} of $X$ given by the \hk metric $g$ is denoted by $TW(X)_\omega$. Its base space is $\mathbb{P}^1$. In the Grassmaniann description of the period domain the periods of a twistor  family are all positive oriented 2 planes contained in the positive space $\langle\,Re(f(\sigma_X)), Im(f(\sigma_X)),f(\omega)\,\rangle$. Any such set of oriented 2 planes contained in a positive 3 space is called \emph{twistor line}. A \emph{generic twistor line} is a twistor line containing the image of the very general point of the moduli space, \ie any point corresponding to a manifold $X$ with no nontrivial divisors. We remark that the period domain is connected by generic twistor lines, cf. \cite[Proposition 3.7]{huy_tor}.\\
Let $L$ be a lattice and let $v\in L$. The divisor of $v$, denoted $div(v)$, is the integer $l\geq0$ satisfying $(v,L)=l\mathbb{Z}$.\\
 For every \hk manifold $X$ we will view $H_2(X,\mathbb{Z})$ as the dual to the lattice $H^2(X,\mathbb{Z})$, which sits naturally inside $H^2(X,\mathbb{Q})$, and let $D^\vee=(1/div(D))D$ for any divisor $D$.
If $S$ is a $K3$ surface and $D$ is a divisor on $S$ we still denote by $D$ the divisor on $S^{[n]}$ defined by the subschemes of $S$ intersecting $D$. We will also denote by $2\delta$ the exceptional divisor of the map $S^{[n]}\rightarrow S^{(n)}$, by $L_n$ the lattice $H^2(S^{[n]},\mathbb{Z})$ and we will call manifolds of \kntiposp all manifolds deformation equivalent to $S^{[n]}$. The moduli spaces and period domains of such manifolds will be denoted respectively $\mathcal{M}_n$ and $\Omega_n$. If $N\subset M$ is a lattice of signature $(3,l)$, we denote with $\mathcal{M}_N$ and $\Omega_N$ the moduli spaces and periods of marked \hk manifolds $(X,f)$ with $H^2(X)\cong M$ such that $f^{-1}(N^\perp)\subset Pic(X)$. In particular the period domain $\Omega_N$, if nonempty, coincides with the period domain for the lattice $N$.

\section{Wall divisors}
The main objects of interest in the present paper are cones related to the geometry of \hk manifolds. Let $X$ be a \hk manifold. The first cone we are interested in is the positive cone, denoted by $\mathcal{C}_X$, which is the connected component containing a \kahl class of the cone of positive classes. Of paramount importance in this paper is the \kahl cone, which is the cone containing all \kahl classes and it is denoted $\mathcal{K}_X$. Then there is the birational \kahl cone, denoted $\mathcal{BK_X}$, which is the union $\cup f^{-1} \mathcal{K}_{X'}$, where $f$ runs through all birational maps between $X$ and any \hk manifold $X'$. If $X$ is projective, then the closure of the birational \kahl cone is just the movable cone, the closure of the cone of divisors whose linear system has no divisorial base components. Finally there is the Mori cone, or the cone of effective curves, which is the closure of the cone in $H_2(X,\mathbb{R})$ generated by classes of effective curves.

Let $X$ be a \hk manifold and let $\omega$ be a \kahl class. Let $\mathcal{F}$ be a cone inside $H_2(X,\mathbb{Z})$ such that $(\omega,r)\geq 0$ for all $r\in \mathcal{F}$. Then the dual cone $\mathcal{F}^\vee$ is the cone 
\begin{equation}
\{D\in H^{1,1}(X,\mathbb{Z}),\,\text{such that } (D,r)> 0\,\,\text{for all}\,r\in\mathcal{F}\,\, \text{and } (D,\omega)\geq 0\}.
\end{equation}
With this notion the \kahl cone of a projective $X$ is just the dual to the Mori cone.\\
We will need some properties of twistor lines and their liftings to moduli spaces, the following is part of the work done by Huybrechts, Markman and Verbitsky in order to prove a global Torelli theorem for \hk manifolds: 
\begin{lem}\label{lem:twistor_lift}\cite[Proposition 5.4]{huy_tor}
Let $X$ be a \hk manifold and let $f\,:\,H^2(X,\mathbb{Z})\,\rightarrow M$ be a marking.  
Then every generic twistor line on $\Omega_M$ lifts to a twistor family in $\mathcal{M}_M^0$.
\begin{proof}
Let $TW\subset \Omega_M$ be a generic twistor line and let $\Delta\subset TW$ be a small open neighbourhood of a distinguished point $0\in \Delta$. Since $\mathcal{P}$ is a local isomorphism, $\Delta$ lifts to $\mathcal{M}^0_M$. 
The line $TW$ is generic, therefore the very general element of $\Delta$ is given by a manifold with zero Picard lattice and with $\mathcal{K}=\mathcal{C}$. Let $X_t$ be one such very general element and let $W$ be the positive 3-space corresponding to $TW$. Let $\alpha_t$ be a class contained in $W\cap\langle Re(\sigma_{X_t}), Im(\sigma_{X_t})\rangle^\perp$. This class is positive and of type $(1,1)$. Therefore $\pm\alpha_t$ lies in the \kahl cone. Moreover $\mathcal{P}(TW(X_t)_{\alpha_t})=TW$.
\end{proof} 
\end{lem}
\begin{defn}\label{defn:muraja}
Let $X$ be a \hk manifold and let $D$ be a divisor on $X$. Then $D$ is called a \emph{\murajasp} if $D^2<0$ and $f^{-1}\circ g(D^\perp)\cap\mathcal{BK}_X=\emptyset$, for all markings $f$ and $g$ of $X$ such that $f^{-1}\circ g$ is a parallel transport Hodge isometry.
\end{defn}
We remind that a parallel transport Hodge isometry is obtained by taking parallel transport of a marking $f$ in a path inside a deformation of $X$. For manifolds of \kntiposp, the group of parallel transport operators has been determined by Markman \cite{mar2}.

The following makes this class of divisors particularly useful.

\begin{thm}\label{thm:muri_kahl}
Let $(X,f)$ and $(Y,g)$ be two marked \hk manifolds and let $D$ be a \murajasp of $X$. Suppose that $(g^{-1}\circ f)(D)\in Pic(Y)$ and that $(X,f)$ and $(Y,g)$ lie in the same connected component $\mathcal{M}^0_n$ of $\mathcal{M}_n$. Then $(g^{-1}\circ f)(D)$ is a \murajasp on $Y$.
\begin{proof}
To prove our claim, we will suppose that there exists a marked manifold $(Z,h)$ such that $D':=h^{-1}\circ f(D)\in Pic(Z)$ and there exists a \kahl class orthogonal to $D'$. We will then deform this manifold with twistor families while keeping $D'$ algebraic. In particular, we will reach a manifold with very general Hodge structure, \ie $Pic=\mathbb{Z}D'$. This is possible since, by contradiction, we have a \kahl class (and therefore an open set of \kahl classes) orthogonal to $D'$. Therefore all the symplectic forms of manifolds belonging to the twistor family associated to this class will be orthogonal to $D'$.\\
Finally, using standard results on connectivity by twistor lines (cf. \cite[Proposition 3.7]{huy_tor} or \cite{beau_tor} for the original reference), we reach a manifold $X'$ birational to $X$, where the algebraic class $D'$ is just $D$. Then we conclude that the manifold $Z$ does not exist, otherwise $X$ would have a class in $\mathcal{BK}_X$ orthogonal to $D$.\\
Let $f(D)=:l$ and let $l^\perp=:N$. Let $\Omega_N$ be the period domain for the lattice $N$ and let $\mathcal{M}^0_N$ be the restriction of $\mathcal{P}^{-1}(\Omega_N)$ to $\mathcal{M}^0_n$. The space  $\mathcal{M}^0_N$ might have several connected components, but if we identify points corresponding to Hodge isometric birational manifolds it is connected. 
The restriction of the period map to $\mathcal{M}_N^0$ is still a local isomorphism onto $\Omega_N$.\\ Assume for contradiction that there exists a point $(Z,h)$ in $\mathcal{M}_N^0$ such that $\mathcal{K}p_Z:=\mathcal{K}_Z\cap h^{-1}(N\otimes\mathbb{R})\neq \emptyset$. Then $\mathcal{K}p_Z$ is an open subset of $h^{-1}(N\otimes\mathbb{R})\cap H^{1,1}(Z)$. Let $\omega_Z$ be a nonrational \kahl class in this intersection. The twistor family $TW(Z)_{\omega_Z}$ is generic, therefore $\mathcal{K}p$ is nonempty also on the very general point of $\mathcal{M}_N^0$.
Let $(Z',h')$ be one such very general point. Then $Pic(Z')=h'^{-1}(l)$ and $\mathcal{K}_{Z'}=\mathcal{C}_{Z'}$, because $h'^{-1}(l)^\perp$ is the only possible wall of the \kahl cone and by our assumption on $Z$ it is not a wall for $Z'$. Then a twistor line passing through $\mathcal{P}(Z',h')\in\Omega_N$ lifts to a generic twistor family inside $\mathcal{M}^0_N$ as in \Ref{lem}{twistor_lift}. Since $\Omega_N$ is connected by generic twistor lines ($N$ has signature (3,19)), we have a chain of twistor families inside $\mathcal{M}_N^0$ connecting $Z$ with a manifold birational to $X$. Therefore $\mathcal{BK}_X\cap f^{-1}(N)\neq \emptyset$. But this is absurd, since $D$ is a \murajasp on $X$. 
\end{proof}
\end{thm}

\begin{lem}\label{lem:estremo_eff}
Let $X$ be a \hk manifold and let $R$ be an extremal ray of the Mori cone of $X$ of negative self intersection. Let $H_2(X,\mathbb{Z})\subset H^2(X,\mathbb{Q})$ by duality. Then any divisor $D$ satisfying $D\in \mathbb{Q}R$ is a \muraja.
\begin{proof}
Since $R$ is an extremal ray, the wall $R^\perp\subset H^2(X,\mathbb{R})$ defines a wall of the \kahl cone of $X$ which is also a wall for the birational \kahl cone. When $X$ is projective, this is can be seen using the flops in the cone theorem \cite[Theorem 3.7]{koll}. If $X$ is not necessarily projective, we can apply \cite[Corollary 5.2]{huy_base} to obtain a cycle $\gamma=Z+\sum_{i}Y_i$ in $X\times X'$ such that $\gamma_{*}^{-1}(\mathcal{K}_{X'})\subset \mathcal{BK}_X$. Here $X'$ is a \hk manifold birational to $X$, $Z$ is the graph of a birational map between the two manifolds and the projections of $Y_i$ are the indeterminacy loci of the map. In our case none of the projections of the $Y_i$ are divisors. If $R$ is not contained in any $Y_i$, then $Z_{*}(R)$ is well defined and effective on $X'$, hence the form $q(\gamma_{*}D,-)$ is positive on the \kahl cone of $X'$. Otherwise, let $\pi(Y_i)$ be the projection to $X$ of $Y_i$. As in \cite[Theorem 7.1]{huy_base}, a \kahl class of $X'$ has negative intersection with a curve in the fibre $Y_i\rightarrow \pi(Y_i)$ and this class is a multiple of $R$. Hence, $q(\gamma_{*}D,-)$ is negative on $\mathcal{K}_{X'}$. \\ Now we need to prove that $h(D)$ is a wall of $\mathcal{BK}_X$ for any parallel transport operator $h$. In the special case of a wall which is also a wall for $\overline{\mathcal{BK}}_X$, our claim is the content of \cite[Section 5]{mar1}. If on the other hand $\mathcal{BK}_X=\mathcal{C}_X$, all Hodge isometries are birational transformations and we are likewise done.\\  
If this is not the case, we can proceed like in \Ref{thm}{muri_kahl} and deny our claim, then we deform with twistor families to a manifold $(Y,g)$ where $\mathcal{BK}_Y=\mathcal{C}_Y$, while keeping the class of the extremal ray algebraic. On $Y$ our claim holds, therefore it does so on $X$.
\end{proof}
\end{lem}
From this we obtain the following:
\begin{prop}\label{prop:mori}
Let $X$ be a \hk manifold and let $\omega$ be a \kahl class. Then there is a bijective correspondence between extremal rays of the Mori cone with negative self intersection and primitive \murajesp $D$ satisfying the additional conditions
\begin{equation}
D^\perp\cap \overline{\mathcal{K}}_X\neq \{0\},\,\,\, (D,\omega)>0\,\,\text{for some}\,\,\omega\in\mathcal{K}_X
\end{equation}
\begin{proof}
Let $R$ be an extremal ray of the Mori cone. Inside $\mathbb{Q}R$ we can choose a primitive divisor $D$ such that $(D,\omega)>0$ by \Ref{lem}{estremo_eff}. Moreover $D^\perp\cap \overline{\mathcal{K}}_X$ consists of all Nef divisors $T$ such that $(T,R)=0$.\\ Conversely, any \murajasp whose perpendicular intersects the closure of the \kahl cone and which is positive on a \kahl class has a dual curve $D^\vee$ which lies in the Mori cone. This ray is on the boundary, since any small modification of $D^\vee$ has orthogonal which meets the interior of the \kahl cone or which no longer meets the closure of the \kahl cone. Being the dual to a wall, this class is also extremal. 
\end{proof}

\end{prop}

Notice that the above proposition exploits some ideas already contained in \cite[Section 5]{huy_ann}, where in particular there is a geometric interpretation of subvarieties containing extremal rays.
This proposition, together with \Ref{thm}{muri_kahl}, provides a numerical characterization of extremal rays of the Mori cone and, dually, walls of the \kahl cone. In the following we analyze some of these rays.\\ 
In particular for manifolds of \kntipo, we will need to analyze the orbit up to isometry of \murajesp inside the lattice $L_n$, but this is made extremely easy by the following, known as Eichler criterion \cite[Lemma 3.5]{ghs1}.

\begin{lem}
\label{lem:ghs_orbit}
Let $T$ be a lattice such that $T\cong U^2\oplus N$ for some lattice $N$ and let $v,w\in T$ be two primitive vectors such that the following conditions hold:
\begin{itemize}
\item $v^2=w^2$.
\item $div(v)= div(w)=m$.
\item $[\frac{v}{m}]=[\frac{w}{m}]$ in $T^\perp/T$.
\end{itemize}
Then there exists an isometry $g$ of $T$ such that $g(v)=w$.
\end{lem}
This lemma gives a unique isometry orbit (if $n-1$ is a prime power) for any primitive element $v$ of $L_n\cong U^3\oplus E_8(-1)^2\oplus (-2n+2)$ as soon as we fix $v^2$ and $div(v)$ and such isometries are of parallel transport.

\section{On manifolds of \kntipo}

\begin{ex}\label{ex:menodue}
Let $S$ be a $K3$ surface and let $C$ be a smooth rational curve on $S$. Then the divisor $C$ on $S^{[n]}$ is a \murajasp and the dual rational curve is an extremal ray of the Mori cone of $X$. Notice that $[C]$ has square $-2$ and divisor $1$.
\end{ex}
\begin{ex}\label{ex:delta}
Let $S$ be a $K3$ surface and let $2\delta$ be the exceptional divisor of the resolution $S^{[n]}\rightarrow S^{(n)}$. Then $\delta$ is a \murajasp and its dual rational curve is an extremal ray of the Mori cone. Notice that $\delta$ has square $-2(n-1)$ and divisor $2(n-1)$. 
\end{ex}
Let $\Lambda=U^4\oplus E_8(-1)^2$ be the Mukai lattice and let $X$ be a manifold of \kntipo. The manifold $X$ comes with a primitive embedding $i\colon H^2(X,\mathbb{Z})\rightarrow \Lambda$, canonical up to a composition with an isometry of $\Lambda$ \cite[Theorem 2.1]{mar2}. A generator $v$ of the rank one lattice orthogonal to the image of $i$ satisfies $v^2=2n-2$. Associated to $D\in H^2(X,\mathbb{Z})$ of square $2-2n$, we get two primitive isotropic classes $w_1$ and $w_2$, which are rational multiples of $v+i(D)$ and $v-i(D)$.
\begin{thm}\cite[Theorem 1.11 and Proposition 1.5]{mar2}\label{thm:mark_muri}
Let $X$ be a manifold of \kntipo. Let $D$ be a divisor on $X$, then the following hold:
\begin{itemize}
\item If $D^2=-2$, then $D$ is a \muraja.
\item If $D^2=2-2n$, $div(D)$ is a multiple of $n-1$ and one of the two isotropic classes $w_1,w_2$, associated to $D$ as above, satisfies $(v,w_i)=1$ or $2$, then $D$ is a \muraja.
\end{itemize}

\end{thm}

Let $M_v(S,H)$ be the moduli space of stable sheaves with primitive Mukai vector $v$ on the $K3$ surface $S$ with respect to a given $v$ generic polarization $H$. Let $s$ be a divisor such that $s^\perp$ gives a wall of the \kahl cone of $M_v(S,H)$ and let $T$ be the rank two hyperbolic primitive sublattice of $\Lambda$ containing $v$ and $s$. Let $\mathcal{P}_T$ be the set $\{t\in T,\,\text{such that } (t,v)>0,\,t^2\geq 0\}$.
\begin{thm}\cite[Theorem 5.7 and 12.1]{bm}\label{thm:eman_muri}
Let $s,v$ and $T$ be as above, then one of the following holds
\begin{itemize}
\item There exists $w\in T$ such that $w^2=-2$ and $(w,v)=0$.
\item There exists $w\in T$ such that $w^2=0$ and $(w,v)=1$ or $2$.
\item There exists $w\in T$ such that $w^2=-2$ and $0<(w,v)\leq v^2/2$.
\item There exist $w,t\in \mathcal{P}_T$ such that $v=w+t$.
\end{itemize}
Moreover for any lattice $T$ satisfying one of the above we obtain a wall of the \kahl cone.
\end{thm} 
Notice that the first two cases are the ones we obtain from \Ref{thm}{mark_muri}. A generalization of the above theorem to arbitrary projective deformations has recently been proved in \cite{bht} by Bayer, Hassett and Tschinkel.
\begin{oss}
Let $T$ be as above and let $D=v^\perp\cap T$ be a wall divisor. Let $d=div(D)$ in $H^2(M_v(S,H))$. A direct computation shows that $d=1$ implies $D^2=-2$ and $d=2$ implies $D^2=-(v^2+8)$ if $n\leq 3$, $D^2\equiv-(v^2+8)$ modulo $4$ otherwise.
\end{oss}

\begin{ex}\label{ex:p2}
Let $S$ be a degree $2$ $K3$ surface with polarization $H$, and let $i$ be the involution of the covering $S\rightarrow \mathbb{P}^2$. Then $\mathbb{P}^2\subset S^{[2]}$, where the projective plane is given by the set of points $(x,i(x))$. Let $l$ be the class of a line in this plane. we have $l=H-3\delta^\vee$, where $2\delta$ is the exceptional divisor. The dual divisor $2H-3\delta$ is a \muraja, since the Mukai flop on the given plane is the reflection along this element.
\end{ex}

\begin{cor}\label{cor:ht_necessario}
Let $X$ be a manifold of \kntiposp and let $R$ be an extremal ray of the Mori cone of $X$. Then $R^2\geq -\frac{n+3}{2} $.
\begin{proof}
Bayer and Macr\`{i} \cite[Proposition 12.6]{bm} proved that the inequality holds on a codimension 2 subset of any connected component of the moduli space of manifolds of \kntipo. This subset is given by moduli space of stable sheaves on projective K3 surfaces. Suppose by contradiction that $R$ is an extremal ray with $R^2<-\frac{n+3}{2}$. 
 Let $D$ be the dual divisor to $R$ (\ie $D/div(D)=R$). Since $R$ is an extremal ray of the Mori cone, we have that $D$ is a \muraja. Let $f$ be a marking of $X$ and let $N=f(D)^\perp$. The space $\mathcal{M}_N$ intersects the codimension 2 subset given by Bayer and Macr\`{i}, therefore $D$ cannot be a \murajasp by \Ref{thm}{muri_kahl}.
\end{proof}
\end{cor}
\begin{ex}\cite[Example 4.11]{ht4}\label{ex:pn}
Let $S$ be a $K3$ surface containing a nonsingular rational curve $E$. Then $E^{[n]}\subset S^{[n]}$ and let $l$ be the class of a line contained in $E^{[n]}$. A direct computation shows $l=E-(n-1)\delta^\vee$. 
 By \Ref{cor}{ht_necessario} $l$ is an extremal ray of the Mori cone, therefore by \Ref{lem}{estremo_eff} the dual divisor $2E-\delta$ is a \murajasp and has square $-(2n+6)$ and divisor $2$. 
\end{ex}
\begin{cor}
Let $X$ be a manifold of \kntiposp where $n-1$ is a prime power. Then divisors of square $-(2n+6)$ and divisor $2$ are \muraje.
\begin{proof}
Let $D$ be one such divisor on $X$ and let $f$ be a marking. The lattice $N:=f(D)^\perp$ is univoquely determined by \Ref{lem}{ghs_orbit} 
 and there is also a single isometry orbit of $f(D)$. In order to obtain our claim using \Ref{thm}{muri_kahl} we need only to find a manifold of \kntiposp where such a divisor is a \muraja, but this is precisely \Ref{ex}{pn}. 
\end{proof}
\end{cor}
\begin{ex}\cite{bm2}
Let $S$ be a $K3$ surface with $Pic(X)=\mathbb{Z}H$, $H^2=2d$ and let $X=S^{[n]}$, $n\geq \frac{d+3}{2}$. Bayer and Macr\`{i} \cite[Proposition 9.3]{bm2} proved that the nef cone of $X$ is generated by $H$ and $H-\frac{2d}{d+n}\delta$. The Mori cone is thus generated by $\delta^\vee=\frac{1}{2(n-1)}\delta$ and $R=H+\frac{d+n}{2n-2}\delta$. The dual divisor $D=H+(d+n)\delta$ is thus a \muraja. Moreover they prove that the curve $R'=R+\delta^\vee$ is not in the Mori cone, even though $R'^{2}\geq -\frac{n+3}{2}$ for large $n$. The least value of $n$ is attained for $d=2$ and $n=5$.  
\end{ex}
\begin{ex}\label{ex:pn-1_bun}
Let $X=M_v(S,H)$ be a moduli space of stable sheaves with $v^2=4$ and such that $X$ contains a primitive divisor $D$ with $D^2=-36$ and $div(D)=4$. We remark that the surjectivity of the period map tells us that the existence of such a $X$ is equivalent to the existence of a primitive embedding $\langle D, H\rangle\hookrightarrow L_n$. If we embed $D$ in the Mukai lattice $\Lambda$  we see that we have a hyperbolic lattice $T$ generated by $v$ and $a=\frac{D+v}{4}$. Since $a^2=-2$ and $(a,v)=1<v^2/2$ we can apply \Ref{thm}{eman_muri} to $T$ and we obtain that $D^\perp$ is a wall of the \kahl cone of $X$. In \cite[Example 4.9]{ht4} such a divisor $D$ is obtained as a $\mathbb{P}^2$ bundle over a $K3$ surface. The same can be done for $n=4,\, v^2=6$ with divisors $D$ satisfying $D^2=-24$ and $div(D)=3$, again represented by $\mathbb{P}^{n-1}$ bundles over $K3$ surfaces.
\end{ex}
\subsection{The case $n=2$}
\begin{prop}
Let $X$ be a manifold of \ktipo. Then the walls of $\mathcal{K}_X$ are given by divisors of square $-2$ or divisors $D$ of square $-10$ and $div(D)=2$.
\begin{proof}
Markman proved that the closure of the birational \kahl cone is dual to the cone generated by divisors of square $-2$, more precisely it is the chamber containing a \kahl class of this cone. 
Hassett and Tschinkel proved that the \kahl cone contains the dual cone to elements of square $-10$ and divisor $2$ and to elements of square $-2$. 
We need only to prove that all divisors with square $-10$ and divisor $2$ give a wall of the \kahl cone (\ie they are \muraja). We remark that such elements are unique up to a parallel transport isometry of $L_2$. \Ref{ex}{p2} tells us that we can apply \Ref{thm}{muri_kahl} to these divisors, therefore we obtain our claim. 
\end{proof}
\end{prop}
Dually, this corollary proves that negative extremal rays of the Mori cone on manifolds of \ktiposp are given by curves of self intersection $-5/2,-2$ or $-1/2$.

\subsection{The case $n=3$}
As an application of \Ref{thm}{muri_kahl} we compute the \kahl and Mori cone of manifolds of \kntiposp for $n=3$. By \Ref{cor}{ht_necessario} the Mori cone is contained in the cone generated by integer classes $r$ in $H_2(X,\mathbb{Z})$ with square greater than $-3$. Looking at their dual divisors not all values between $-3$ and $0$ can be reached, for example a divisor $D$ with $div(D)=2$ has $D^2\equiv -4$ modulo $8$ and if $div(D)=4$ we have $D^2\equiv -4$ modulo $32$. A direct computation shows that the only possible values are those in the following table.
\begin{table}[ht]\label{tab:raggi3}
\begin{tabular}{|c|c|c|}
\hline
$r^2$ & $D^2$ & $div(D)$\\
\hline
$-3$ & $-12$ & $2$\\
\hline
$-9/4$ & $-36$ & $4$\\ 

\hline
$-2$ & $-2$ & $1$\\ 
\hline
$-1/4$ & $-4$ & $4$\\ 
\hline
$-1$ & $-4$ & $2$\\ 
\hline
\end{tabular}
\end{table}
The first four elements in this table are those contained in \cite[Table H3]{ht4}. The last element is one of the possibilities in \Ref{thm}{mark_muri}.\\  
\begin{defn}
Let $X$ be a manifold of $K3^{[3]}$ type and let $\omega$ be a \kahl class. Let $\mathcal{NK}_X\subset\mathcal{C}_X$ be the dual cone to the closure of the cone generated by divisors $D$ having the same numerical invariants as those in \Ref{tab}{raggi3} and satisfying $(D,\omega)\geq 0$.
\end{defn}
\begin{thm}\label{thm:cono3}
Let $X$ be a manifold of $K3^{[3]}$ type. Then $\mathcal{NK}_X= \mathcal{K}_X$.
\begin{proof}
The cone $\mathcal{NK}_X$ is contained in the \kahl cone because all of the walls of $\mathcal{K}_X$ are given by a divisor in the above table. Conversely the cone $\mathcal{NK}_X$ contains the \kahl cone because the elements of the above table are \murajesp by \Ref{ex}{pn-1_bun}, \Ref{ex}{pn}, \Ref{ex}{delta}, \Ref{ex}{menodue} and \Ref{thm}{mark_muri}.   
\end{proof}
\end{thm}
\subsection{The case $n=4$}
As a further application of \Ref{thm}{muri_kahl}, we compute the \kahl and Mori cone of manifolds of \kntiposp for $n=4$. By \Ref{cor}{ht_necessario} the Mori cone is contained in the cone generated by integer classes $r$ in $H_2(X,\mathbb{Z})$ with square greater than $-\frac{7}{2}$. As in the previous case not all values can be obtained, below is a list of all possible extremal rays and of their corresponding dual divisors.
\begin{table}[ht]\label{tab:raggi4}
\begin{tabular}{|c|c|c|}
\hline
$r^2$ & $D^2$ & $div(D)$\\
\hline
$-7/2$ & $-14$ & $2$\\
\hline
$-8/3$ & $-24$ & $3$\\
\hline
$-13/6$ & $-78$ & $6$\\
\hline
$-2$ & $-2$ & $1$\\
\hline
$-2/3$ & $-6$ & $3$\\
\hline
$-1/6$ & $-6$ & $6$\\
\hline
$-3/2$ & $-6$ & $2$\\
\hline
\end{tabular}
\end{table}
The first six elements in this table are those contained in \cite[Table H4]{ht4}. 
Let us remark that this is the highest dimension where the conjecture of Hassett and Tschinkel holds.\\
As before let $X$ be a manifold of $K3^{[4]}$ type with a \kahl class $\omega$. Let $\mathcal{NK}_X$ be the dual cone to the cone generated by divisors $D$ with the same numerical invariants as the elements in the above table and satisfying $(D,\omega)\geq 0$.
\begin{thm}
Let $X$ be a manifold of $K3^{[4]}$ type. Then $\mathcal{NK}_X= \mathcal{K}_X$.
\begin{proof}
The proof goes exactly as in \Ref{thm}{cono3}, the only thing we need to prove is that elements of square $-6$ and divisor $2$ are \muraje.\\
Let $X=M_v(S,H)$ be a smooth moduli space of stable sheaves on a $K3$ surface $S$ with $v^2=6$ and let $D\in Pic(X)$ be an element satisfying $D^2=-6$ and $div(D)=2$. The existence of such a manifold is again due to the surjectivity of the period map. If we embed $D$ in $\Lambda$ we obtain a hyperbolic lattice $T$ generated by $v$ and $a=\frac{v+D}{2}$. The lattice $T$ satisfies the conditions of \Ref{thm}{eman_muri}, therefore by \Ref{thm}{muri_kahl} elements in its numerical equivalence class are \muraje.
\end{proof}
\end{thm}

\bibliographystyle{amsplain}

\end{document}